\documentclass[10pt]{article}
\usepackage[latin1]{inputenc}

\usepackage{amsmath, amssymb, amsthm, amsfonts, amsgen}

\usepackage{indentfirst}

\usepackage{graphicx}

\numberwithin{equation}{section}

\makeatletter

\def\rightharpoonupfill@{\arrowfill@\relbar\relbar\rightharpoonup}

\newcommand{\xrightharpoonup}[2][]{\ext@arrow 0359\rightharpoonupfill@{#1}{#2}}

\def\debaixodolim#1#2{\mathrel{}\mathop{\lim}\limits^{#1}_{#2}}

\makeatother

\newcommand{\ds}{\displaystyle}

\newcommand{\Nb}{{\mathbb{N}}}

\newcommand{\Qb}{{\mathbb{Q}}}

\newcommand{\Rb}{{\mathbb{R}}}

\newcommand{\Zb}{{\mathbb{Z}}}

\def\leq{\leqslant}

\def\geq{\geqslant}

\let\e=\varepsilon

\let\O=\Omega

\let\o=\omega

\let\G=\Gamma

\let\go=\rightarrow

\newtheorem{thm}{Theorem}[section]

\newtheorem{defi}[thm]{Definition}

\newtheorem{rmk}[thm]{Remark}

\newtheorem{lemma}[thm]{Lemma}

\newtheorem{proposition}[thm]{Proposition}

\newtheorem{coro}[thm]{Corollary}

\begin{document}

\title{\bf 3D-2D analysis of a thin film with periodic microstructure}

\author{ Jean-Fran\c cois Babadjian and Margarida Ba\'{\i}a}

\date{}

\maketitle

\begin{abstract}

\noindent The  purpose of this article  is to study  the behavior
of a heterogeneous thin film whose microstructure oscillates on a
scale that is comparable to that of the thickness of the domain.
The argument is based on a 3D-2D dimensional reduction through a
$\Gamma$-convergence analysis, techniques of two-scale convergence
and a decoupling procedure between the oscillating variable and the in-plane variable.\\

\noindent {\bf Keywords:}  dimension reduction, thin films,
periodic  integrands, $\G$-convergence, two-scale convergence,  quasiconvexity, equi-integrability.\\

\noindent {\bf MSC 2000 (${\cal A}{\cal M}{\cal S}$):} 35E99,
35M10, 49J45, 74B20, 74K35, 74G65, 74Q05

\end{abstract}

\section{Introduction and main result}

\noindent In this work we study the asymptotic behavior of a
heterogeneous $\varepsilon$-thin domain with periodic
microstructure of period $\e$, as  $\e\go 0$, through a
$\Gamma$-limit analysis.
 Techniques of two-scale convergence and a decoupling procedure
between the microscopic oscillating variables and the macroscopic
in-plane variables are used to derive the relaxed two-dimensional
energy from its three-dimensional counterpart.

Let $\o$ be an open and bounded subset of $\Rb^{2}$. For each
$0<\e \ll 1$ define  $\Omega_{\e}:=\omega\times (-\e,\e)$.
Consider a deformable thin body occupied by a hyperelastic
material with a periodic microstructure of period $\e$ whose
reference configuration is given by the thin domain
$\Omega_{\varepsilon}$, and whose stored energy density
$W(\e):\O_\e \times \Rb^{3 \times 3} \to \Rb$  is assumed to be  a
Carath\'eodory function satisfying some $p$-growth and coercivity
conditions ($1<p<\infty$). We assume that the body is pinned on
the lateral boundary, that is $v(x)=x \; \text{ on }
\partial \omega \times (-\e,\e)$, for all its admissible
deformations, and that it is  subjected to the action of regular
surface traction densities $g(\e)$ on $\Sigma_\e:=\omega \times
\{-\e,\e\}$, and regular dead loads $f(\e)$. The  total energy of
this body under the action of this forces is the difference
between the elastic energy and the work of external forces.  More
precisely,

$$\mathcal E(\e)(v):= \int_{\O_\e}W(\e)(x;Dv)\, dx - \int_{\Omega_\e}f(\e) \cdot v\, dx
-\int_{\Sigma_\e}g(\e) \cdot v \, d \mathcal H^2, $$

\noindent for  $v\in \mathcal V(\e):= \{v \in W^{1,p}(\O_\e;\Rb^3)
: v(x)=x \; \text{ on } \partial \omega \times (-\e,\e)\}$, and
where $\mathcal H^2$ stands for the two-dimensional Hausdorff
measure. It may occur that the minimization problem associated
with this energy admits no solution over the set of kinematically
admissible fields $\mathcal V(\e)$. However, we can introduce the
notion of \textit{almost-minimizer} of $\mathcal E(\e)$, $v(\e)
\in \mathcal V(\e)$, by requiring that

$$\mathcal E(\e)(v(\e)) \leq \inf_{v \in \mathcal V(\e)} \mathcal E(\e)(v)+\e \,
h(\e),$$

\noindent where $h(\e) \searrow 0^+$ when $\e \to 0$. Note that if
the minimization problem admits a solution -- for instance if
$W(\e)$ is quasiconvex in its second variable -- then we can take
$h \equiv 0$.

As usual, in order to study this problem as $\e\go 0$ we rescale
the $\varepsilon$-thin body into a reference domain of unit
thickness (see e.g.\! Acerbi, Buttazzo and Percivale \cite{ABP},
Anzellotti, Baldo and Percivale \cite{AnBP}, Le Dret and Raoult
\cite{LDR},
 Braides,  Fonseca and Francfort \cite{Bra&Fo&Fr}), so that the resulting energy will be
defined on a fixed body, while the dependence on $\varepsilon$
turns out to be explicit in the transverse derivative. For this,
we consider the change of variables

$$ \Omega_{\varepsilon}\rightarrow \Omega:=\omega\times I, \quad (x_1, x_2, x_3)
\mapsto  \left(x_1, x_2, \frac{1}{\varepsilon}x_3\right),$$

\noindent  and define $u(x_\alpha,x_3 / \e)=v(x_\alpha,x_3)$ on
the rescaled cylinder $\O$, where $I:=(-1,1)$  and
$x_\alpha:=(x_1,x_2)$ is the in-plane variable. It is well known
 that membrane theory
arises at the order $\e$ of a formal asymptotic expansion (see
Fox, Raoult and Simo \cite{FRS}), provided that the body forces
are of order $1$ and the surface loadings are of order $\e$. Since
this energy is of order $\e$ we divide the total energy by $\e$
and, in addition we  assume that

$$\left\{
\begin{array}{rcl}
f(\e)(x_\alpha,\e x_3) & = & f(x_\alpha,x_3),\\
g(\e)(x_\alpha,\e x_3) & = & \e \, g(x_\alpha,x_3),
\end{array}
\right.$$

\noindent where $f \in L^{p'}(\Omega;\Rb^3)$, $g \in
L^{p'}(\Sigma;\Rb^3)$ ($1/p+1/p'=1$) and $\Sigma:=\o \times \{
-1,1\}$.  If $W_\e(x_\alpha,x_3;\cdot)=W(\e)(x_\alpha,\e
x_3;\cdot)$, for fixed $\e$ minimizing $\mathcal E(\e)$ on
$\mathcal V(\e)$ is equivalent to minimizing

$$\mathcal E_\e(u):=\frac{\mathcal E(\e)(v)}{\e}=\int_\O W_\e
\left(x;D_\alpha u(x)\Big| \frac{1}{\e}D_3 u(x)\right)\, dx -
\int_\O f \cdot u\, dx -\int_\Sigma g \cdot u \, d\mathcal H^2$$

\noindent on $\mathcal V_\e:=\{u \in W^{1,p}(\O;\Rb^3) :
u(x)=(x_\alpha,\e x_3) \text{ on }\partial \omega \times I\}$.
Denote by $D_i=\frac{\partial}{\partial x_i}$ for $i \in
\{1,2,3\}$ and $D_\alpha= (D_1,D_2)$.   In the sequel, we identify
$\Rb^{d\times N}$ (resp. $\Qb^{d \times N}$) with the space of
real (resp. rational) $d\times N$ matrices. For all
$\overline{\xi}=(z_1|z_2)\in \Rb^{3\times 2}$ and $z\in \Rb^3$,
$(\overline{\xi}|z)$ is the matrix whose first two columns are
$z_1$ and $z_2$ and whose last one is $z$. Denoting a
almost-minimizer of the  rescaled energy  by
$u_\e(x_\alpha,x_3):=v(\e)(x_\alpha,\e x_3)$, we obtain

\begin{equation}\label{quasimin}
\mathcal E_\e (u_\e) \leq \inf_{u \in \mathcal V_\e} \mathcal
E_\e(u) +h(\e).
\end{equation}

Our aim is to study the asymptotic behavior of the equilibrium
problem (\ref{quasimin}) as $\e\go 0$ via a $\G$-convergence
method (we refer to Braides and Defranceschi \cite{Bra&De},
Braides \cite{B} and  Dal Maso \cite{DM} for a comprehensive
treatment and bibliography on $\G$-convergence).

The motivation for studying problem  (\ref{quasimin}) comes from
the work in  Braides, Fonseca and Francfort \cite{Bra&Fo&Fr} where
the authors have established an abstract dimensional reduction
variational convergence result in a general setting  for a family
of stored energies of the form $W_{\e}(x;\xi)$ and derived
specific characterizations for particular cases. In Section 3 of
\cite{Bra&Fo&Fr}  a heterogeneous nonlinear
 membrane model is derived by $\Gamma$-convergence, and  heterogeneity in the
transverse
 direction is considered. Precisely, the authors treat the case where  the stored
energy density is of the form
  $ W(x_{3};\xi)$,
generalizing the previous work of Le Dret and Raoult in \cite{LDR}
who treated a homogeneous material,  i.e.\!  when $W$ depends only
in $\xi$. Later, Babadjian and Francfort \cite{Ba&Fr} considered
energies of the form $ W(x;\xi)$ with a general heterogeneity.
Furthermore  in Section 4 of  \cite{Bra&Fo&Fr}, a 3D-2D analysis
coupled with a homogenization in the in-plane direction is studied
in the case where $W_{\e}(x;\xi)= W(x_{3}, x_{\alpha}/\e;\xi) $.
Shu \cite{S} also investigated similar problems, in the framework
of martensitic materials, with different length scales for the
film thickness and the material microstructure.

Here we propose  to establish a dimensional reduction and
homogenization result, where both scales are identical, by adding
in the stored energy density an explicit dependence on the
macroscopic in-plane variable $x_\alpha$. Namely, we assume that
$W_\e(x_\alpha,x_3;\cdot)=W(x_\alpha,x_3,x_\alpha /\e;\cdot)$ for
some function $W: \O\times \Rb^{2}\times \Rb^{3\times
3}\rightarrow \Rb$ whose hypotheses will be introduced later.

 Two features differentiate our approach from what is available in most of the
literature in the subject. The first one is the use of a two-scale
convergence argument (see Nguetseng \cite{N1,N2} and Allaire
\cite{A} for the notion and properties of two-scale convergence).
The same argument was used by Ba\'{\i}a and Fonseca in \cite{B&Fo}
in a purely  homogeneous context, i.e. without considering the
dimensional reduction problem.   The second feature is due to the
definition of the homogenized stored energy in which two
independent variables  appear  (see identity (\ref{whom}) below).
To take into account this structure, we are led to decouple the
macroscopic in-plane variable $x_{\alpha}$ from  the microscopic
oscillating variable $x_{\alpha}/ \e$  via  an extension
 argument along  the lines of Babadjian and
Francfort \cite{Ba&Fr}.

For a comprehensive treatment on the homogenization of integral
functionals  via a $\G$-limit approach, we refer to Braides and
Defranceschi \cite{Bra&De} and references therein. We will denote
by $\mathcal L^{N}$ the $N$-dimensional Lebesgue measure in
$\Rb^{N}$\! (in the sequel $N$ will be equal to $2$ or $3$).

For each $\e>0$ we define $\mathcal I_\varepsilon
:L^p(\Omega;\mathbb R^3)\rightarrow \overline{\Rb}$ by

\begin{equation}\label{defIe}\mathcal I_\varepsilon(u):=\left\{\begin{array}{ll}
\ds\int_{\O}W\left(x_{\alpha},x_{3},\frac{x_\alpha}{\varepsilon};D_\alpha
u(x) \Big| \frac{1}{\varepsilon}D_3 u(x)\right) dx
 & \text{if } u\in
W^{1,p}(\O;\mathbb R^3),\\\\
+\infty &\text{otherwise},
\end{array}\right.\end{equation}

\noindent with $1< p<\infty$, where we assume that   $W:\O \times
\mathbb R^2 \times \mathbb R^{3 \times 3} \rightarrow \Rb$
satisfies the following hypotheses:

\begin{itemize}

\item[$(H_{1})$] $W(x,\,\cdot\,;\,\cdot\,)$ is continuous for
a.e.\! $x\in \O$;

 \item[$(H_{2})$] $W(\,\cdot\,,y_{\alpha};\xi)$ is  measurable for all $(y_{\alpha},\xi) \in \Rb^{2}\times \mathbb R^{3 \times 3}$;

\item[$(H_{3})$] there exists $0<\beta <+ \infty$  such that

$$\frac{1}{\beta}|\xi|^p -\beta \leq W(x,y_{\alpha};\xi) \leq \beta(1+|\xi|^p),
\quad {\rm for\,\, a.e.\,\,} x \in \Omega\, \text{ and\,
for\,all\,}\,  (y_{\alpha},\xi) \in \Rb^{2} \times \mathbb R^{3
\times 3};$$

\item[$(H_{4})$] $W(x,\,\cdot\,;\xi)$ is $Q'$-periodic for a.e.\!
$x \in \Omega$ and all $\xi \in \mathbb R^{3 \times 3}$, where we
denote by $Q'=(0,1)^2$ the unit cube of $\mathbb R^2$\!.
\end{itemize}

\begin{rmk}\label{1}

{\rm  We remark that due to hypothese $(H_{1})$ and $(H_{2})$ the
function $W$ is a  Carath\'{e}odory integrand as
$W(x,\cdot;\cdot)$ is continuous a.e.\! $x\in \O$ and
$W(\cdot,y_{\alpha};\xi)$ is measurable for all $y_{\alpha}\in
\Rb^{2}$ and $\xi\in \Rb^{3\times 3}$. This implies (see e.g.\!
Proposition 3.3 in Braides and Defranceschi \cite{Bra&De} or
Proposition 1.1, Chapter VIII in Ekeland and Temam \cite{Ek&Te})
that $W$ is equivalent to a Borel function, that is there exist a
Borel function $\tilde W$ such that
$W(x,\,\cdot\,;\,\cdot\,)=\tilde W(x,\,\cdot\,;\,\cdot\,)$ for
a.e.\! $x\in \O$. As a consequence the integral in (\ref{defIe})
is well defined.
 As noted by Allaire in \cite{A}, Section 5, the
measurability of $W$ in the pair $(x,y_{\alpha})$ does not let us
conclude that, for fixed $\xi$, the function $x \mapsto
W\left(x,x_{ \alpha}/\e;\xi \right)$ is measurable. The continuity
of $W(x,y_\alpha;\xi)$ in at least one of the variables $x$ or
$y_\alpha$ turns out to be sufficient to guarantee the
measurability of this function. In the present paper, we decide to
impose the continuity in the $y_\alpha$ variable. Note that we
could also have considered $W$ to be continuous in $x$ and
measurable in $y_{\alpha}$ but the proof of our main result does
not hold anymore in this context. }
\end{rmk}

As for notation, we will identify $W^{1,p}(\o;\Rb^{3})$ with the
set of functions $u\in W^{1,p}(\O;\Rb^{3})$ such that
$D_{3}u(x)=0$ for a.e.\!  $x \in \O$ and we will use the  notation
$\G (L^{p}(\O))$-limit whenever we refer to the $\G$-convergence
with respect to the usual metric in $L^{p}(\O;\Rb^3)$.  We prove
the following main result.

\begin{thm}\label{jf-m} If \,$W$ satisfies $(H_1)$-$(H_{4})$,
then the family $\{\mathcal I_\e\}_{\e >0}$
$\G(L^p(\O))$-converges to  the functional $\mathcal I_{\rm hom}:
L^{p}(\Omega;\mathbb R^3)\rightarrow \overline\Rb$ defined by

\begin{equation}
{\mathcal I}_{\rm hom} (u):= \left\{\begin{array}{ll} \ds
2\int_{\o} W_{\rm hom}(x_{\alpha};D_{\alpha} u(x_{\alpha}))\,
 dx_{\alpha} & \mbox{if}\,\, u\in W^{1,p}(\o;\Rb^{3}), \\&\\
+ \infty & \mbox{otherwise},
\end{array}\right.
\end{equation}

where $W_{\rm hom}$ is given by

\begin{eqnarray}\label{whom}\displaystyle
W_{\rm hom}(x_\alpha;\overline \xi) & \hspace{-1cm}:=
\hspace{-1cm}& \lim_{T \to +\infty} \hspace{-0.1cm}\inf_\varphi
\Bigg\{ \frac{1}{2 T^2}  \int_{(0,T)^2 \times
I}\hspace{-0.5cm}W\big(x_\alpha,y_3,y_{\alpha};\overline \xi
+D_\alpha \varphi(y) |D_3 \varphi(y)\big)
dy \hspace{-0.1cm}: \nonumber \\
& & \varphi \in W^{1,p}((0,T)^{2} \times I;\mathbb R^3), \;
\varphi=0 \text{ on }\partial (0,T)^{2} \times I \Bigg\}
\end{eqnarray}

\noindent for  a.e.\! $x_\alpha \in \omega$ and all $\overline \xi
\in \mathbb R^{3 \times 2}$.
\end{thm}

As a consequence of  Theorem \ref{jf-m}  we deduce the usual
convergence of (almost-)minimizers. More precisely, we have the
following result.

\begin{coro}

Let $\{u_\e\} \subset \mathcal V_\e$ be a sequence of
almost-minimizers for $\{\mathcal I_{\e}\}_{\e>0}$ (see identity
(\ref{quasimin})). Then $\{u_\e\}$ is weakly relatively compact in
$W^{1,p}(\O;\Rb^3)$. Furthermore, any limit point $u$ of this
sequence is a solution of the minimization problem

$$\min_{v-(x_\alpha,0) \in W^{1,p}_0(\omega;\Rb^3)}\left\{ 2 \int_\o W_{\rm
hom}(x_\alpha;D_\alpha
 v(x_\alpha))\, dx_\alpha - \int_\o(\overline f + g^+ + g^-)(x_\alpha)\cdot
v(x_\alpha)\, dx_\alpha\right\},$$ where $\overline f :=
\frac{1}{2}\int_I f(\cdot,x_3)\, dx_3$ and $g^\pm:=g(\cdot,\pm
1)$.

\end{coro}

This corollary departs from the classical result  on the type of
boundary conditions that have been considered (see e.g.\
Proposition 7.2 in Braides and Defranceschi  \cite{Bra&De}). This
difficulty is overcome by the fact that  we can prescribe the
lateral boundary conditions of recovering sequences (see Remark
\ref{bci}).  We do not include the proof of this corollary here
because it is similar to that of Corollary 1.3 in Bouchitt\'e,
Fonseca and Mascarenhas \cite{BFM2}.

The plan of this work is as follows: In Section 2 we will discuss
some properties of $W_{\rm hom}$, namely that it is well defined,
proving that the limit on the right hand side of (\ref{whom})
exists, and that $W_{\rm hom}(x_{\alpha};\,\cdot\,)$  is
continuous for a.e. $x_{\alpha}\in \o$. Section 3 is devoted to
the proof of our main result, Theorem \ref{jf-m}. The starting
point of our analysis is the $\Gamma$-limit integral
representation result, Theorem 2.5,  in Braides, Fonseca and
Francfort \cite{Bra&Fo&Fr}. Our objective is to identify the limit
integrand, showing that it coincides (almost everywhere) with
$W_{\rm hom}$. We will use an argument of two-scale convergence to
derive an upper bound for the limit integrand (Lemma \ref{ineq1}).
Since the problem at fixed $\e$ and the asymptotic problem  as
$\e\go 0$ are of different nature (one is a full three-dimensional
problem, the other a two-dimensional one),
 we will need
to use a decoupling argument to prove the other  inequality (Lemma
\ref{ineq2}).  For this purpose
 it will be convenient to  extend $W$ to a function which  is (separately) continuous
everywhere. This is the aim of  Lemma \ref{ext} (see Appendix in
Section 4) which provides conditions under which a Carath\'eodory
function such as $W$ can be extended to a separately continuous
function in the macroscopic in-plane variable $x_{\alpha}$ and the
microscopic variable $x_{\alpha}/\e.$

\section{Preliminary results}

\noindent In this section we will prove some properties of  the
stored energy $W_{\rm hom}$ that will be of use in the proof of
Theorem \ref{jf-m}.

\begin{rmk}\label{Wnonneg}{\rm

To prove Theorem \ref{jf-m} we may assume, without loss of
generality, that $W$ is non negative. Indeed, in view of $(H_{3})$
it suffices to replace $W$ by $W+\beta$. }\end{rmk}

We begin by showing that in the definition (\ref{whom}) of $W_{\rm
hom}$  the limit as $T \to +\infty$ exists. The proof of this
property is a direct consequence of a result due to Licht and
Michaille \cite{LM}, Theorem 3.1 (see also Lemma 4.3.6 in
Bouchitt\'{e}, Fonseca and  Mascarenhas  \cite{BFM}).

\begin{lemma}\label{welldef}

If\,  $W$ satisfies  $(H_1)$-$(H_4)$, then

\begin{eqnarray*}
W_{\rm hom}(x_\alpha;\overline \xi) \hspace{-0.3cm}& =
\hspace{-0.3cm}& \lim_{T \to +\infty} \inf_{\varphi} \Big\{
\frac{1}{2 T^2} \int_{(0,T)^2 \times
I}\hspace{-0.3cm}W\big(x_\alpha,y_3,y_{\alpha};\overline \xi
+D_\alpha \varphi(y) |D_3 \varphi(y)\big)
dy: \\
& & \hspace{2cm}\varphi \in W^{1,p}((0,T)^{2} \times I;\mathbb
R^3), \; \varphi=0 \text{ on }\partial (0,T)^{2} \times I \Big\}
\end{eqnarray*}

\noindent  exists for a.e.\! $x_{\alpha}\in \o$ and all
$\overline{\xi}\in \Rb^{3\times 2}$\!.

\end{lemma}

\noindent {\it Proof.} Let $x_{\alpha}\in \o$ be such that
$(H_1)$, $(H_{3})$ and $(H_{4})$  hold and let $\overline{\xi}\in
\Rb^{3\times 2}$\!. Define $\mu: {\cal A}(\Rb^{2})\go \Rb^{+}$ by

\begin{eqnarray*}
\mu(A) & := & \inf_\varphi \Big\{ \frac{1}{2}\int_{A \times
I}W(x_\alpha,y_3,y_{\alpha};\overline \xi +D_\alpha \varphi(y)|D_3
\varphi(y))\,
dy: \\
& & \hspace{1cm} \varphi \in W^{1,p}(A\times I;\mathbb R^3), \;
\varphi=0 \text{ on }\partial A \times I \Big\},
\end{eqnarray*}

\noindent  where $\mathcal A(\Rb^2)$ stands for the family of open
subsets of $\Rb^2$.

The function $\mu$ is well defined and, thanks to $(H_{3})$, it is
finite. Moreover this set function satisfies the assumptions of
Theorem 3.1 in Licht and Michaille \cite{LM}. Indeed firstly, by
$(H_{3})$, $\mu(A)\leq \beta(1+|\overline \xi|^p)\,{\cal
L}^{2}(A)$ for all $A\in {\cal A}(\Rb^{2})$. Secondly, $\mu$ is
subadditive, that is $\mu(C)\leq \mu(A)+\mu(B)$ for all $A,\, B,\,
C \in {\cal A}(\Rb^{2})$ with $A \cap B \neq \emptyset$ and
$\overline{C}=\overline{A}\cup\overline{B}$. Finally, by
$(H_{4})$, for any $\mathbf i \in \Zb^2$, $\mu(A+\mathbf
i)=\mu(A)$ for all $A\in {\cal A}(\Rb^{2})$. As a consequence  the
limit

$$\lim_{T\go +\infty} \frac{\mu((0,T)^{2})}{T^{2}}=W_{\rm
hom}(x_{\alpha};\overline{\xi}) $$

\noindent exists. \hfill$\blacksquare$

\vspace{0.2cm}

\begin{rmk}\label{infwhom}
{\rm It can be proved that the limit as $T \to +\infty$ in
(\ref{whom}) can
 be replaced by an infimum taken for every $T>0$ (see Braides and Defranceschi
\cite{Bra&De} or Ba\'{\i}a and Fonseca \cite{B&Fo}). }
\end{rmk}

Now that $W_{\rm hom}$ is well defined, we will show that
 $W_{\rm hom}(x_\alpha; \cdot)$ is continuous for a.e.\! $x_\alpha
\in \o$ for later use in Theorem \ref{jf-m}. To prove this
property directly it seems that  we  need  more than merely the
continuity condition imposed on $W(x,y_\alpha;\cdot)$ (e.g.\! a
$p$-Lipschitz condition). We remark that if $W(x,y_\alpha;\cdot)$
was  quasiconvex, then by the $p$-growth condition $(H_3)$,
$W(x,y_\alpha;\cdot)$  would satisfy a $p$-Lipschitz condition
(see Lemma \ref{cont} below). Since  we do not want to a priori
restrict  the stored energy density too much, in order to
compensate for this lack of regularity we  first prove in Lemma
\ref{QWW} that the value of $W_{\rm hom}$ does not change if we
replace $W$ by its quasiconvexification ${\mathcal Q}W$ (see
Remark \ref{qsgrmk} below).

\begin{rmk}\label{qsgrmk}
{\rm For a.e.\! $x\in \O$, all $y_{\alpha}\in \Rb^{2}$ and all
$\xi\in\Rb^{3\times 3}$ define

$$
{\mathcal Q}W(x,y_{\alpha};\xi):=[
QW(x,y_{\alpha};\,\cdot\,)](\xi)$$

\noindent where $ QW(x,y_{\alpha};\,\cdot\,)$ stands for the usual
quasiconvexification of $W(x,y_{\alpha};\,\cdot\,)$. Then, the
function $\mathcal QW(x,y_{\alpha};\,\cdot\,)$ is quasiconvex (see
e.g.\! Dacorogna \cite{Dac}) and  satisfies $(H_1)$-$(H_4)$  with
the exception that $\mathcal QW(x,\cdot;\xi)$ may only be  upper
semicontinuous  for a.e.\! $x \in \O$ and all $\xi \in \Rb^{3
\times 3}$\! (as the infimum of continuous functions).  By an
argument similar to that of Lemma \ref{welldef} we conclude that

\begin{eqnarray*}
(\hspace{-0.05cm}{\mathcal Q}W\hspace{-0.05cm})_{\rm
hom}(x_\alpha;\overline \xi) \hspace{-0.3cm}& = &
\hspace{-0.4cm}\lim_{T \to +\infty} \hspace{-0.1cm}\inf_{\varphi}
\Bigg\{ \hspace{-0.05cm}\frac{1}{2 T^2}\hspace{-0.05cm}
\int_{(0,T)^2 \times I}\hspace{-0.5cm}{\mathcal
Q}W\big(x_\alpha,y_3,y_{\alpha};\overline \xi
\hspace{-0.05cm}+\hspace{-0.05cm}D_\alpha \varphi(y) |D_3
\varphi(y)\big)\,
dy \hspace{-0.05cm}: \\
& & \hspace{1.5cm}\varphi \in W^{1,p}((0,T)^{2} \times I;\mathbb
R^3), \; \varphi=0 \text{ on }\partial (0,T)^{2} \times I \Bigg\}
\end{eqnarray*}

  \noindent  exists for a.e.\! $x_\alpha \in \o$ and all $\overline{\xi}\in \Rb^{3\times
2}$\!. }

\end{rmk}

\begin{lemma}\label{QWW}

If $W$ satisfies $(H_1)$-$(H_4)$, then $({\mathcal Q}W)_{\rm hom}
(x_\alpha;\overline \xi)=W_{\rm hom} (x_\alpha;\overline \xi)$ for
a.e.\! $x_{\alpha}\in \o$ and all $\overline \xi\in \Rb^{3\times
2}.$

\end{lemma}

\noindent {\it Proof. } Let $x_\alpha \in \omega$ be such that
both $({\mathcal Q}W)_{\rm hom}(x_\alpha;\,\cdot\,)$ and $W_{\rm
hom}(x_\alpha;\,\cdot\,)$ are well defined. Since $W \geq
{\mathcal Q}W$, we  have $W_{\rm hom}(x_\alpha;\overline \xi)\geq
({\mathcal Q}W)_{\rm hom}(x_\alpha;\overline \xi)$  for all
$\overline \xi\in \Rb^{3\times 2}$\!. Let us prove now the
converse inequality. Let $\overline\xi\in  \Rb^{3\times 2}$\!. For
each $n>0$, let $T_{n} \in \mathbb N$ and $\varphi_{n} \in
W^{1,\infty}((0,T_n)^2 \times I;\mathbb R^3)$ satisfying
$\varphi_{n}=0$ on $\partial (0,T_{n})^2 \times I$, be such that

$$({\mathcal Q}W)_{\rm hom}(x_\alpha;\overline \xi) +
\frac{1}{n}\geq \frac{1}{2 T_{n}^2} \int_{(0,T_{n})^2 \times
I}{\mathcal Q}W(x_\alpha,y_3,y_\alpha;\overline \xi +D_\alpha
\varphi_{n}(y)|D_3 \varphi_{n}(y))\, dy.$$

The Lipschitz regularity of $\varphi_{n}$ is ensured  because of
the density of $W^{1,\infty}((0,T_n)^2 \times I;\mathbb R^3)$ in
$W^{1,p}((0,T_n)^2 \times I;\mathbb R^3)$ together with the
$p$-growth condition $(H_{3})$. Thus

\begin{equation}\label{1707}({\mathcal Q}W)_{\rm hom}(x_\alpha;\overline \xi)\geq
\limsup_{n \to +\infty} \frac{1}{2 T_{n}^2} \int_{(0,T_{n})^2
\times I}{\mathcal Q}W(x_\alpha,y_3,y_\alpha;\overline \xi
+D_\alpha \varphi_{n}(y)|D_3 \varphi_{n}(y))\, dy.\end{equation}

 For each $n \in \Nb$ fixed, by Acerbi-Fusco Relaxation Theorem (see Lemma III.1 and
Statement III.7 in \cite{AF}) and  Remark \ref{Wnonneg}, there
exists a sequence $\{\varphi_{n,k}\}_{k}\subset
W^{1,\infty}((0,T_{n})^2 \times I;\Rb^{3})$  satisfying
$\varphi_{n,k}=\varphi_n$ on $\partial [(0,T_n)^2 \times I]$ with
$\varphi_{n,k}\xrightharpoonup[k\go +\infty]{} \varphi_{n}$  and
such that

$$\begin{array}{ll}&\ds \frac{1}{2 T_{n}^2} \int_{(0,T_{n})^2 \times I}{\mathcal
Q}W(x_\alpha,y_3,y_\alpha;\overline \xi +D_\alpha
\varphi_{n}(y)|D_3 \varphi_{n}(y))\, dy\\\\
&\ds =\quad \lim_{k \to +\infty}\frac{1}{2 T_{n}^2}
\int_{(0,T_{n})^2 \times I}W(x_\alpha,y_3,y_\alpha;\overline \xi
+D_\alpha \varphi_{n,k}(y)|D_3 \varphi_{n,k}(y))\,
dy.\end{array}$$

 \noindent From (\ref{1707}) we have

\begin{eqnarray*}
({\mathcal Q}W)_{\rm hom}(x_\alpha;\overline
\xi)\hspace{-0.2cm}&\geq &\hspace{-0.2cm}\limsup_{n \to
+\infty}\limsup_{k \to +\infty}\frac{1}{2 T_{n}^2}
\int_{(0,T_{n})^2 \times
I}\hspace{-0.4cm}W(x_\alpha,y_3,y_\alpha;\overline
 \xi +D_\alpha \varphi_{n,k}(y)|D_3 \varphi_{n,k}(y))\, dy\\
 \hspace{-0.2cm}& \geq &\hspace{-0.2cm} \limsup_{n\go\infty}\,\, \inf_\varphi \Bigg\{ \frac{1}{2 T^2_{n}}
\int_{(0,T_{n})^2 \times
I}\hspace{-0.4cm}W(x_\alpha,y_3,y_{\alpha};\overline \xi +D_\alpha
\varphi(y)|D_3 \varphi(y))\, dy
:\nonumber\\
&&\hspace{2.5cm}\varphi \in W^{1,p}((0,T_{n})^{2} \times I;\mathbb
R^3), \;
\varphi=0 \text{ on }\partial (0,T_{n})^{2} \times I \Bigg\}\nonumber\\
&=& \hspace{-0.2cm}W_{\rm hom}(x_{\alpha};\overline \xi).
\end{eqnarray*}

\hfill$\blacksquare$\\

We are now in position to prove the continuity of $W_{\rm hom}$ in
its second variable :

\begin{lemma}\label{cont}Let $W$ satisfying $(H_1)$-$(H_4)$, then
$W_{\rm hom}(x_{\alpha};\,\cdot\,)$ is continuous on $\Rb^{3\times
2}$ for a.e.\! $x_{\alpha}\in \o$.

\end{lemma}

 \noindent {\it Proof.}  We observe that by
the $p$-growth condition in $(H_3)$ and Remark \ref{qsgrmk},
${\mathcal Q}W$ satisfies a $p$-Lipschitz condition (see
Marcellini \cite{Ma1}): There exists $\beta>0$ such that for all\
$y_\alpha \in \Rb^2$ and a.e.\! $x \in \O$,

\begin{equation}\label{plipschitz}
|{\mathcal Q}W(x,y_\alpha;\xi_1)-{\mathcal Q}W(x,y_\alpha;\xi_2)|
\leq \beta( 1+|\xi_1|^{p-1}+|\xi_2|^{p-1})|\xi_1-\xi_2|, \quad
\xi_1, \, \xi_2 \in \Rb^{3 \times 3}.
\end{equation}

  Take $x_{\alpha} \in \o$ such that both
$({\mathcal Q}W)_{\rm hom}(x_{\alpha};\,\cdot\,)$ and $W_{\rm
hom}(x_{\alpha};\,\cdot\,)$ are well defined. By  Lemma \ref{QWW}
we have  $({\mathcal Q}W)_{\rm hom}(x_{\alpha};\,\cdot\,)= W_{\rm
hom}(x_{\alpha};\,\cdot\,)$. Given $\overline{\xi}\in \Rb^{3\times
2}$ let $\overline{\xi}_{n}\rightarrow\overline{\xi}$ in
$\Rb^{3\times 2}$\!. From the definition of $W_{\rm
hom}(x_{\alpha};\overline\xi)$,  for fixed $\delta>0$ choose $T
\in \Nb$ and $\varphi \in W^{1,p}((0,T)^{2} \times I;\mathbb R^3),
\; \varphi=0 \text{ on }\partial (0,T)^{2} \times I$, such that

\begin{equation}\label{in}
W_{\rm hom}(x_{\alpha};\overline{\xi})+\delta
 \geq \frac{1}{2T^{2}}\int_{(0,T)^{2}\times I}W(x_\alpha,y_3,y_{\alpha};\overline \xi
+D_\alpha \varphi(y)|D_3 \varphi(y))\,
 dy.\end{equation}

\noindent  Therefore, Remark \ref{infwhom} yields

\begin{eqnarray*}
\limsup_{n \go +\infty} W_{\rm hom}(x_{\alpha};\overline \xi_n) &
\leq & \limsup_{n \go
+\infty}\frac{1}{2T^{2}}\int_{(0,T)^{2}\times
I}W(x_\alpha,y_3,y_{\alpha}; \overline \xi_{n} +D_\alpha
\varphi(y)|D_3
\varphi(y))\,dy\\
& = & \frac{1}{2T^{2}}\int_{(0,T)^{2}\times
I}W(x_\alpha,y_3,y_{\alpha};\overline \xi +D_\alpha \varphi(y)|D_3
\varphi(y))\, dy
\end{eqnarray*}

\noindent due to hypothesis $(H_1),$ the $p$-growth condition in
$(H_3)$ and Lebesgue's Dominated Convergence Theorem. So by
(\ref{in}) and   letting $\delta\go 0$ we conclude that

\begin{equation}\label{upper}\limsup_{n \go +\infty}W_{\rm
hom}(x_{\alpha};\overline{\xi_{n}})\leq W_{\rm
hom}(x_{\alpha};\overline{\xi}).\end{equation} Similarly,  for
each  $n\in \Nb$ consider $T_{n}\in \Nb$ ($T_{n}\nearrow +\infty$)
and $\varphi_{n} \in W^{1,p}((0,T_{n})^{2} \times I;\mathbb R^3),
\; \varphi_{n}=0 \text{ on }\partial (0,T_{n})^{2} \times I$, such
that

\begin{eqnarray*}
W_{\rm hom}(x_{\alpha};\overline{\xi}_{n})+\frac{1}{n}& \geq &
\frac{1}{2T^{2}_{n}}\int_{(0,T_{n})^{2}\times I}{\mathcal
Q}W(x_\alpha,y_3,y_{\alpha};\overline \xi_{n} +D_\alpha
\varphi_{n}(y)|D_3 \varphi_{n}(y))\,
dy\\
&=&\frac{1}{2}\int_{Q'\times I}{\mathcal
Q}W(x_\alpha,y_3,T_{n}y_{\alpha};\overline \xi_{n} +D_\alpha
\varphi_{n}(T_{n}y_\alpha,y_3)|D_3
\varphi_{n}(T_{n}y_\alpha,y_3))\,
dy\\
&=& \frac{1}{2}\int_{Q'\times I}{\mathcal
Q}W(x_\alpha,y_3,T_{n}y_{\alpha};\overline \xi_{n} +D_\alpha
\psi_{n}(y)|T_n D_3 \psi_{n}(y))\, dy,
\end{eqnarray*}

\noindent  after a change of variables and where
$\psi_{n}(y):=\frac{1}{T_{n}}\varphi_{n}(T_{n}y_\alpha,y_3)$.
Clearly the function $\psi_{n}$ belongs to $W^{1,p}(Q'\times
I;\mathbb R^3)$ and $\psi_{n}=0 \text{ on }\partial Q'\times I$.

By the $p$-coercivity hypothesis in $(H_3)$  and (\ref{upper}),
the sequence $\{(D_\alpha \psi_{n}|T_n D_3 \psi_n)\}$ is bounded
in $L^{p}(Q'\times I;\Rb^{3\times 3})$ uniformly in $n$. We can
write that

\begin{eqnarray*}
&&\liminf_{n \go +\infty}\int_{Q' \times I}{\mathcal
Q}W(x_\alpha,y_3,T_{n}y_{\alpha};\overline \xi_{n} +D_\alpha
\psi_{n}(y)|T_{n}D_3 \psi_{n}(y))\,
dy\\
&& \geq  \liminf_{n \go +\infty}\int_{Q' \times
I}\left[\,{\mathcal Q}W(x_\alpha,y_3,T_{n}y_{\alpha};\overline
\xi_{n} +D_\alpha \psi_{n}(y)|T_{n}D_3
\psi_{n}(y))\right.\\
&& \left.\hspace{4.0cm} -{\mathcal
Q}W(x_\alpha,y_3,T_{n}y_{\alpha};\overline \xi
+D_\alpha \psi_{n}(y)|T_{n}D_3 \psi_{n}(y))\,\right]\,dy\\
&&\hspace{0.4cm}+ \liminf_{n \go +\infty}\int_{Q' \times
I}{\mathcal Q}W(x_\alpha,y_3,T_{n}y_{\alpha};\overline \xi
+D_\alpha \psi_{n}(y)|T_{n}D_3 \psi_{n}(y))\, dy.
\end{eqnarray*}

\noindent Using (\ref{plipschitz}), H\"{o}lder inequality, the
fact that $\{\|(D_\alpha \psi_{n}|T_n D_3 \psi_n)\|_{L^{p}(Q'
\times I;\Rb^{3 \times 3})}\}$ is bounded and
$\overline{\xi}_{n}\rightarrow \overline{\xi}$, we obtain

\begin{eqnarray*}
&&\liminf_{n \go +\infty}\int_{Q' \times I}\left[\,{\mathcal
Q}W(x_\alpha,y_3,T_{n}y_{\alpha};\overline \xi_{n} +D_\alpha
\psi_{n}(y)|T_{n}D_3
\psi_{n}(y))\right.\\
&&\hspace{2.0cm}\left.-{\mathcal
Q}W(x_\alpha,y_3,T_{n}y_{\alpha};\overline \xi +D_\alpha
\psi_{n}(y)|T_{n}D_3 \psi_{n}(y))\,\right]\,dy=0,
\end{eqnarray*}

\noindent and consequently

\begin{eqnarray}\label{lower}
\liminf_{n \go +\infty} W_{\rm hom}(x_{\alpha};\overline{\xi}_{n})
& \geq & \liminf_{n \go +\infty}\frac{1}{2}\int_{Q' \times
I}{\mathcal Q}W(x_\alpha,y_3,T_{n}y_{\alpha};\overline \xi
+D_\alpha \psi_{n}(y)|T_{n}D_3 \psi_{n}(y))\,
dy\nonumber\\
& = & \liminf_{n \go
+\infty}\frac{1}{2T^{2}_{n}}\int_{(0,T_{n})^{2}\times I}{\mathcal
Q}W(x_\alpha,y_3,y_{\alpha};\overline \xi +D_\alpha
\varphi_{n}(y)|D_3 \varphi_{n}(y))\,
dy\nonumber\\
&\geq &  ({\mathcal Q}W)_{\rm hom}(x_{\alpha};\overline \xi)\nonumber \\
&=& W_{\rm hom}(x_{\alpha};\overline \xi).
\end{eqnarray}

\noindent From (\ref{upper}) and (\ref{lower}), we conclude that
$W_{\rm hom}(x_{\alpha};\cdot)$ is continuous at $\overline{\xi}$.
\hfill$\blacksquare$

\section{Proof of Theorem \ref{jf-m}}

 \noindent We start by localizing our functionals. Representing by ${\cal A}(\o)$ the
class of all open subsets of $\o$, define $\mathcal I_\varepsilon
:L^p(\Omega;\mathbb R^3) \times \mathcal A(\omega) \rightarrow
\overline\Rb$ by
$$\mathcal I_\varepsilon(u;A):=\left\{\begin{array}{ll} \ds
\int_{A\times
I}W\left(x_{\alpha},x_3,\frac{x_\alpha}{\varepsilon};D_\alpha u(x)
\Big| \frac{1}{\varepsilon}D_3 u(x)\right) dx & \text{if } u\in
W^{1,p}(A \times I;\mathbb R^3),\\&\\
+\infty &\text{otherwise}.
\end{array}\right.$$

We will prove  that the family of functionals  $\{\mathcal
I_\varepsilon(\cdot;A)\}_{\e>0}$ $\G$-converges with respect to
the $L^{p}(A\times I;\Rb^{3})$-topology to the functional
${\mathcal I}_{\rm hom} (\cdot ;A): L^p(\O;\mathbb R^3)
\rightarrow \overline\Rb$

\begin{equation}\label{claim} {\mathcal I}_{\rm hom} (u ;A):=
\left\{\begin{array}{ll} \ds 2\int_{A} W_{\rm
hom}(x_{\alpha};D_{\alpha} u(x_{\alpha}))\,
 dx_{\alpha} \quad & \mbox{if } u\in W^{1,p}(A;\Rb^{3}),\\&\\
+ \infty & \mbox{otherwise},
\end{array}\right.\end{equation}

\noindent for all $A\in {\cal A}(\o)$. As a consequence, taking
$A=\o$ yields Theorem \ref{jf-m}.

For any $A\in  \mathcal A(\omega)$ and any sequence $\{\e_j\}
\searrow 0^+$, consider the $\G$-lower limit of the family
$\{\mathcal I_{\varepsilon_{j}}(\cdot;A)\}_{j \in \Nb}$, $\mathcal
I_{\{\varepsilon_j\}}(\cdot;A): L^p(\O;\mathbb R^3) \rightarrow
\overline\Rb$,  given  by

\begin{equation}\label{F}
\mathcal I_{\{\varepsilon_j\}}(u ;A):=\inf_{\{u_j\}}\left\{
\liminf_{j \rightarrow +\infty} \mathcal I_{\varepsilon_j}(u_j;A)
:  u_j \rightarrow u \text{ in }  L^p(A \times I;\mathbb{R}^3)
\right\}.
\end{equation}

\begin{rmk}\label{dem2}
{\rm In view of  the coercivity condition $(H_{4})$, for all $A\in
{\cal A}(\o)$ we have that $\mathcal
I_{\{\varepsilon_j\}}(u;A)=+\infty\!$ whenever $u \in
L^p(\Omega;\mathbb{R}^3) \setminus W^{1,p}(A;\mathbb{R}^3)$, hence
our objective is to characterize $\mathcal
I_{\{\varepsilon_j\}}(u;A)$ for $u\in W^{1,p}(A;\mathbb R^3)$. }
\end{rmk}

By virtue of Remark \ref{dem2}, together with Theorem 2.5 in
Braides, Fonseca and Francfort  \cite{Bra&Fo&Fr}, it follows that
every sequence $\{\varepsilon_j\}$ admits a subsequence
$\{\varepsilon_{j_n}\} \equiv \{\e_n\}$ such that $\mathcal
I_{\{\varepsilon_n\}}(\,\cdot\,;A)$ defined in (\ref{F}) is the
$\Gamma(L^p(A \times I))$-limit of $\{\mathcal
I_{\varepsilon_n}(\,\cdot\,;A)\}_{n \in \Nb}$ for all $A \in
\mathcal{A}(\omega)$. Further there exists a Carath\'eodory
function $W_{\{\varepsilon_n\}}:\omega \times \mathbb{R}^{3 \times
2} \rightarrow \mathbb{R}$ such that

\begin{equation}\label{intrep}
\mathcal I_{\{\varepsilon_n\}}(u;A)=2\int_A
W_{\{\varepsilon_n\}}(x_{\alpha};D_{\alpha} u(x_{\alpha}))\,
dx_{\alpha},
\end{equation}

\noindent for all $A \in \mathcal{A}(\omega)$ and all $u \in
W^{1,p}(A;\mathbb{R}^3)$.  Our aim is to show that $\mathcal
I_{\{\e_n\}}(\cdot;A)=\mathcal I_{\rm hom}(\cdot;A)$ on
$W^{1,p}(A;\Rb^3)$ for all $A\in {\cal A}(\o)$. Given $A\in {\cal
A}(\o)$, in view of the integral representation (\ref{intrep}) and
(\ref{claim}), it is enough to show that
$W_{\{\e_n\}}(x_\alpha;\overline \xi)=W_{\rm
hom}(x_\alpha;\overline \xi)$ for a.e.\! $x_\alpha \in A$ and all
$\overline \xi\in \Rb^{3\times 2}$, and thus to work with affine
functions instead of general Sobolev functions. We will prove that
$W_{\{\e_n\}}(x_\alpha;\overline \xi)=W_{\rm
hom}(x_\alpha;\overline \xi)$ for a.e.\! $x_\alpha \in \o$ and all
$\overline \xi\in \Rb^{3 \times 2}$.

\begin{rmk}
{\rm \label{bci} Lemma 2.6 of Braides, Fonseca and Francfort
\cite{Bra&Fo&Fr} implies that $\mathcal
I_{\{\varepsilon_j\}}(u;A)$ is unchanged if the approximating
sequences $\{u_j\}$ are constrained to match the lateral boundary
condition of their target, i.e. $u_j\equiv u \mbox{ on }\partial A
\times I$. }
\end{rmk}

>From now onward,  $\{\varepsilon_n\}$ will denote a subsequence of
$\{\varepsilon_j\}$ for which the $\Gamma(L^p(A \times I))$-limit
of $\{\mathcal I_{\varepsilon_n}(\cdot ;A)\}_{n \in \Nb}$ exists
and coincides with $\mathcal I_{\{\varepsilon_n\}}(\cdot ;A)$ for all $A\in {\cal A}(\o)$.\\

For each  $T>0$ consider ${\cal S}_{T}$ a countable set of
functions in $\mathcal C^\infty([0,T]^2 \times [-1,1];\Rb^{3})$
that is dense in

$${\cal W}_T=\{ \varphi \in W^{1,p}((0,T)^{2}\times I;\Rb^{3}):\,\, \varphi=0 \,\text{
on } \,  \partial{(0,T)^{2}} \times I\}.$$

\begin{defi}\label{setL}
{\rm Let $L$  be the set of Lebesgue points $x_{\alpha}^{0}$ for
all functions

\begin{equation}\label{C1)}
W_{\{\e_n\}}(\cdot;\overline{\xi}), \quad W_{\rm
hom}(\cdot;\overline{\xi})
\end{equation}

\noindent and

\begin{equation}\label{C2)}
x_{\alpha}\mapsto \int_{Q' \times I}W(x_{\alpha},y_3,T
y_{\alpha};\overline{\xi}+D_{\alpha}\varphi(Ty_\alpha,y_3)
|D_{3}\varphi(Ty_\alpha,y_3))\, dy,
\end{equation}

\noindent with $T \in \Nb$, $\varphi\in {\cal S}_{T}$ and
$\overline{\xi}\in \Qb^{3\times 2}$, and for which  $W_{\rm
hom}(x_{\alpha}^{0};\,\cdot\,)$ is well defined.}
\end{defi}

We have that $\mathcal L^2(\o\setminus L)=0$. Given
$x_{\alpha}^{0}\in L$, we denote by $Q'(x_\alpha^0,\delta)$ the
cube in $\Rb^2$ centered in $x_\alpha^0$ and of side length
$\delta>0$ where $\delta$ is small enough so that
$Q'(x_\alpha^0,\delta) \in \mathcal A(\o)$.

To prove that $W_{\{\e_n\}}(x_\alpha;\overline \xi)=W_{\rm
hom}(x_\alpha;\overline \xi)$ for a.e.\! $x_\alpha \in \o$ and all
$\overline \xi\in \Rb^{3 \times 2}$ we first show in Lemmas
\ref{ineq1} and \ref{ineq2} below that both functions coincide on
$L \times \mathbb Q^{3 \times 2}$\!. The general case will only be
treated at the end of that section  using the Carath\'eodory
property of both integrands.

Fix $\overline{\xi}\in \Qb^{3\times 2}$ and set $v(x):= \overline
\xi \cdot x_\alpha$. By (\ref{intrep}) and (\ref{C1)})

\begin{eqnarray}\label{lpg}
W_{\{\e_n\}}(x_{\alpha}^{0};\overline{\xi})& = &  \lim_{\delta\go
0}\frac{1}{\delta^{2}}\int_{Q^{\prime}(x_{\alpha}^{0},\delta)}W_{\{\e_n\}}(x_{\alpha};
\overline{\xi})\,dx_{\alpha}\nonumber\\
& = &  \lim_{\delta\go 0}\frac{{\cal
I}_{\{\e_n\}}(v;Q^{\prime}(x_{\alpha}^{0},\delta))}{2\delta^{2}}.
\end{eqnarray}

\begin{lemma}\label{ineq1}
$W_{\{\e_n\}}(x_{\alpha}^{0};\overline\xi)\leq W_{\rm
hom}(x_{\alpha}^{0};\overline\xi)$ for all $x_{\alpha}^{0}\in L$
and all $\overline\xi\in \Qb^{3\times 2}.$
\end{lemma}

\noindent {\it Proof.}   Given $k\in \Nb$, let $T_{k}\in \Nb$ and
${{\varphi}}_k \in {\cal S}_{T_{k}}$ with $\varphi_{k}=0$ on
$\partial(0,T_k)^{2}\times I$, be such that

$$W_{\rm hom}(x_{\alpha}^{0};\overline\xi)+\frac{1}{k} \geq
\frac{1}{2{T_{k}}^{2}}\int_{(0,T_{k})^{2} \times
I}W(x_{\alpha}^{0},y_3,y_{\alpha}; \overline\xi+D_{\alpha}
{{\varphi}}_k(y)| D_{3} {{\varphi}}_k(y))\, dy.$$ \noindent This
is possible because of the continuity properties $(H_{1})$ of $W$,
the growth conditions $(H_{3})$ and the density of ${\cal
S}_{T_{k}}$ in ${\cal W}_{T_{k}}$.  Extend $\varphi_{k}$
periodically with period $T_{k}$ to $\Rb^{2}\times I$.  For $x\in
\Rb^{2} \times I$, define $u^{k}_{n}(x):=\overline\xi\cdot
x_{\alpha}+\e_{n}\varphi_{k}(\frac{x_{\alpha}}{\e_{n} },x_{3})$.

For fixed $k$, $u^{k}_{{n}} \rightarrow v$ in
$L^{p}(Q^{\prime}(x_{\alpha}^{0},\delta)\times I;\Rb^{3})$ as
$n\go \infty$, hence, by (\ref{lpg})

\begin{eqnarray*}
\ds W_{\{\e_n\}}(x_{\alpha}^{0};\overline\xi) \hspace{-0.3cm}&
\leq &\hspace{-0.3cm} \liminf_{\delta \go 0}\liminf_{n \go
+\infty}
\frac{1}{2\delta^{2}}\int_{Q^{\prime}(x_{\alpha}^{0},\delta)\times
I}\hspace{-0.5cm}
W\left(\hspace{-0.05cm}x_{\alpha},x_3,\frac{x_{\alpha}}{\e_{n}};D_\alpha
u_n^k(x) \Big |
\frac{1}{\varepsilon_n}D_3 u_n^k(x)\hspace{-0.05cm}\right) dx\\
\hspace{-0.3cm}&\ds = &\hspace{-0.3cm}\liminf_{\delta \go
0}\liminf_{n \go \hspace{-0.05cm}+\infty}\hspace{-0.05cm}
\frac{1}{2\delta^{2}}\hspace{-0.05cm}\int_{Q'(x_\alpha^0,\delta)
\times
I}\hspace{-0.5cm}W\hspace{-0.05cm}\left(\hspace{-0.05cm}x_{\alpha},x_3,\hspace{-0.05cm}\frac{x_\alpha}{\varepsilon_n};\overline
\xi \hspace{-0.05cm}+\hspace{-0.05cm}D_\alpha
\varphi_k\hspace{-0.06cm}\left(\hspace{-0.05cm}\frac{x_{\alpha}}{\e_{n}},x_{3}\hspace{-0.05cm}\right)\hspace{-0.1cm}\Big|D_3
\varphi_k\hspace{-0.06cm}\left(\hspace{-0.05cm}\frac{x_{\alpha}}{\e_{n}},x_{3}\hspace{-0.05cm}\right)\hspace{-0.1cm}\right)\hspace{-0.05cm}
dx.
\end{eqnarray*}

\noindent Define

$$h_{k}(x_\alpha,y_\alpha):=\int_{-1}^{1}W(x_{\alpha},x_3,
T_{k}y_{\alpha};\overline\xi+D_{\alpha}\varphi_{k}(T_{k}y_{\alpha},x_{3})|D_{3}\varphi_
{k} (T_{k}y_{\alpha},x_{3}))dx_3, $$

 \noindent for a.e.\,  $x_\alpha \in \o$ and  $y_\alpha \in
\Rb^{2}$.  The continuity of $W$ with respect to $y_{\alpha}$, its
measurability and periodicity  properties, and the fact that
$T_{k}\in \Nb$ lead us to conclude that   the function $h_{k}\in
L^{1}(Q'(x_{\alpha}^{0},\delta); \mathcal C_{\rm per}(Q'))$ for
fixed $\delta>0$, where $\mathcal C_{\rm per}(Q')$ denotes the
space of $Q'$-periodic and continuous functions defined on
$\Rb^2$\! (see Lemma 5.3 in Allaire \cite{A}). Lemma 5.2 in
\cite{A} together with Fubini's Theorem  yields to

\begin{eqnarray}\label{1320}
&& \lim_{n \go +\infty}\int_{Q'(x_\alpha^0,\delta) \times
I}W\left(x_{\alpha},x_3,\frac{x_\alpha}{\varepsilon_n};\overline
\xi +D_\alpha
\varphi_k\left(\frac{x_{\alpha}}{\e_{n}},x_{3}\right)\Big|D_3
\varphi_k\left(\frac{x_{\alpha}}{\e_{n}},x_{3}\right)\right)\,
dx\nonumber\\
& = &\hspace{-0.2cm} \lim_{n \go +
\infty}\int_{Q'(x_\alpha^0,\delta)}
h_{k}\left(x_\alpha,\frac{x_{\alpha}}{T_{k}\e_{n}}\right)\, dx_\alpha\nonumber\\
& =  &\hspace{-0.2cm}\int_{Q'(x_{\alpha}^{0},\delta)}\int_{Q'}
h_{k}(x_\alpha,y_\alpha)\, dy_\alpha\,
dx_\alpha\nonumber\\
& = &\hspace{-0.2cm} \int_{Q'(x_{\alpha}^{0},\delta)}\int_{Q'
\times I}\hspace{-0.2cm} W(x_{\alpha},x_3, T_{k}y_{\alpha};
\overline\xi+D_{\alpha}\varphi_{k}(T_{k}y_{\alpha},x_{3})|D_{3}\varphi_{k}
(T_{k}y_{\alpha},x_{3})) dy_{\alpha}\, dx_{3}\, dx_{\alpha}.
\end{eqnarray}

\noindent Using  (\ref{C2)}) and passing to the limit in
(\ref{1320}), as $\delta\to 0$, we have that

\begin{eqnarray*}
&& W_{\{\e_n\}}(x_{\alpha}^{0};\overline\xi)\\
&&  \leq  \frac{1}{2 } \int_{Q' \times I} W(x_{\alpha}^0,x_3, T_k
y_{\alpha}; \overline\xi+D_{\alpha}\varphi_{k}(T_k
y_{\alpha},x_{3})|D_{3}\varphi_{k}
(T_k y_{\alpha},x_{3})) dy_{\alpha}\,dx_{3}\\
&&  \leq  W_{\rm hom}(x_{\alpha}^{0};\overline\xi)+\frac{1}{k}.
\end{eqnarray*}

\noindent

Letting $k\go +\infty$ we assert the claim.

\hfill$\blacksquare$\\

Note that the same proof could be used to prove Lemma 2.5 in
Babadjian and Francfort \cite{Ba&Fr}.

\begin{lemma}\label{ineq2}

$W_{\{\e_n\}}(x_{\alpha}^{0};\overline\xi)\geq W_{\rm
hom}(x_{\alpha}^{0};\overline\xi)$ for all $x_{\alpha}^{0}\in L$
and all $\overline\xi\in \Qb^{3\times 2}.$

\end{lemma}

\noindent {\it Proof.}  Let $\{v_n\} \subset
W^{1,p}(Q'(x_{\alpha}^{0},\delta) \times I;\Rb^3)$ be a recovering
sequence for the $\Gamma$-limit, i.e.

$$v_n \to 0 \text{ in }L^p(Q'(x_{\alpha}^{0},\delta) \times I;\Rb^3)$$

and

$$\mathcal I_{\{\varepsilon_n\}}(v;Q'(x_\alpha^0,\delta)) = \lim_{n \to +\infty}
\int_{Q'(x_\alpha^0,\delta) \times
I}W\left(x_\alpha,x_3,\frac{x_\alpha}{\varepsilon_n};\overline \xi
+D_\alpha v_n(x) \Big| \frac{1}{\varepsilon_n} D_3
v_n(x)\right)dx.$$ According to Theorem 1.1 in Bocea and Fonseca
\cite{Bo&Fo}, there exists a subsequence of $\{\varepsilon_n\}$
(not relabelled) and a sequence $\{u_n\} \subset
W^{1,p}(Q'(x_\alpha^0,\delta) \times I;\mathbb R^3)$ such that,
upon setting $E_n:=\{x \in Q'(x_\alpha^0,\delta) \times I :
u_n(x)=v_n(x)\}$, we have that

\begin{equation}\label{bf}
\left\{ \begin{array}{l}
u_n \to 0 \text{ in }L^p(Q'(x_\alpha^0,\delta) \times I;\mathbb R^3),\\\\
\left\{ \left| \left( D_\alpha u_n \big| \frac{1}{\varepsilon_n}
D_3 u_n\right)
\right|^p \right\} \text{ is equi-integrable},\\\\
\ds \lim_{n \to +\infty} \mathcal L^3([Q'(x_\alpha^0,\delta)
\times I]\setminus E_n)=0.
\end{array} \right.
\end{equation}

Thus, in view of the $p$-growth condition $(H_{3})$ together with
(\ref{bf}) and Remark \ref{Wnonneg},

\begin{eqnarray*}
\mathcal I_{\{\varepsilon_n\}}(v;Q'(x_\alpha^0,\delta)) & \geq &
\limsup_{n \to +\infty}
\int_{E_n}W\left(x_\alpha,x_3,\frac{x_\alpha}{\varepsilon_n};\overline
\xi
+D_\alpha u_n \Big| \frac{1}{\varepsilon_n} D_3 u_n\right)dx\\
& = & \limsup_{n \to +\infty} \int_{Q'(x_\alpha^0,\delta) \times
I}W\left(x_\alpha,x_3,\frac{x_\alpha}{\varepsilon_n};\overline \xi
+D_\alpha u_n \Big|
\frac{1}{\varepsilon_n} D_3 u_n\right)dx\\
& & -\limsup_{n \to +\infty} \int_{[Q'(x_\alpha^0,\delta) \times
I]\setminus
E_n}W\left(x_\alpha,x_3,\frac{x_\alpha}{\varepsilon_n};\overline
\xi +D_\alpha u_n
\Big| \frac{1}{\varepsilon_n} D_3 u_n\right)dx\\
& \geq & \limsup_{n \to +\infty} \int_{Q'(x_\alpha^0,\delta)
\times
I}W\left(x_\alpha,x_3,\frac{x_\alpha}{\varepsilon_n};\overline \xi
+D_\alpha u_n \Big| \frac{1}{\varepsilon_n} D_3 u_n\right)dx.
\end{eqnarray*}

For any $h \in \mathbb N$, we split $Q'(x_\alpha^0,\delta)$ into
$h^2$ disjoints cubes $Q'_{i,h}$ of side length $\delta/h$ so that
$Q'(x_\alpha^0,\delta)=\bigcup_{i=1}^{h^2} Q'_{i,h}$ and

\begin{equation}\label{equiint}
\mathcal I_{\{\varepsilon_n\}}(v;Q'(x_\alpha^0,\delta)) \geq
\limsup_{h \to +\infty} \limsup_{n \to +\infty}
\sum_{i=1}^{h^2}\int_{Q'_{i,h} \times
I}W\left(x_\alpha,x_3,\frac{x_\alpha}{\varepsilon_n};\overline \xi
+D_\alpha u_n \Big| \frac{1}{\varepsilon_n} D_3 u_n\right)dx.
\end{equation}

For every $\eta>0$ and $\lambda>0$, let $K_\eta \subset \Omega$
and $W^{\eta,\lambda}$ be given by Lemma \ref{ext} below (with
$N=d=3$, $m=2$ and $f=W$). Then

\begin{equation}\label{mes}
\mathcal L^3(\Omega \setminus K_\eta) < \eta.
\end{equation}

On the other hand, define

$$R^\lambda_n:=\left\{ x \in Q'(x_\alpha^0,\delta) \times I: \left|\left(\overline \xi
+D_\alpha u_n(x) \Big| \frac{1}{\varepsilon_n} D_3
u_n(x)\right)\right| \leq \lambda\right\}.$$

Chebyshev's inequality implies that there exists a constant $C>0$
-- which does not depend on $n$ or $\lambda$ -- such that

\begin{equation}\label{lambda}
\mathcal L^3 ([Q'(x_\alpha^0,\delta) \times I] \setminus
R^\lambda_n) < \frac{C}{\lambda^p}.
\end{equation}

In what follows  we denote by $\ds \limsup_{\lambda,\eta,h,n}$ the
successive $\ds \limsup_{\lambda \to +\infty} \limsup_{\eta \to 0}
\limsup_{h \to +\infty} \limsup_{n \to +\infty}$. Since $W$ and
$W^{\eta,\lambda}$ coincide on $K_\eta \times \Rb^2 \times
\overline{B}(0,\lambda)$,  where in the sequel the set
$\overline{B}(0,\lambda)$ stands for the closed ball $\{ \xi \in
\Rb^{3 \times 3}: |\xi|\leq \lambda\}$ of $\Rb^{3 \times 3}$, we
get in view of (\ref{equiint})

\begin{eqnarray*}
&\ds \mathcal I_{\{\varepsilon_n\}}(v;Q'(x_\alpha^0,\delta)) \geq\\
&\ds \limsup_{\lambda, \eta, h, n} \sum_{i=1}^{h^2}
\int_{[Q'_{i,h} \times I] \cap R^\lambda_n \cap
K_\eta}W^{\eta,\lambda}\left(x_\alpha,x_3,\frac{x_\alpha}{\varepsilon_n};\overline
\xi +D_\alpha u_n \Big| \frac{1}{\varepsilon_n} D_3 u_n\right)dx.
\end{eqnarray*}

By virtue of  (\ref{ContExt}) below  and (\ref{mes}),

$$\sum_{i=1}^{h^2} \int_{([Q'_{i,h} \times I] \cap R^\lambda_n )\setminus
K_\eta}W^{\eta,\lambda}\left(x_\alpha,x_3,\frac{x_\alpha}{\varepsilon_n};\overline
\xi +D_\alpha u_n \Big| \frac{1}{\varepsilon_n} D_3 u_n\right) dx
\leq \beta(1+\lambda^p)\eta \xrightarrow[\eta \to 0]{}0,$$

uniformly in $(n,h)$, so that

\begin{eqnarray*}
&\ds \mathcal I_{\{\varepsilon_n\}}(v;Q'(x_\alpha^0,\delta)) \geq \\
&\ds \limsup_{\lambda,\eta,h,n} \sum_{i=1}^{h^2} \int_{[Q'_{i,h}
\times I] \cap
R^\lambda_n}W^{\eta,\lambda}\left(x_\alpha,x_3,\frac{x_\alpha}{\varepsilon_n};\overline
\xi +D_\alpha u_n \Big| \frac{1}{\varepsilon_n} D_3 u_n\right)dx.
\end{eqnarray*}

\noindent Fix $y_\alpha \in Q'$.  Since
$W^{\eta,\lambda}(\cdot,y_\alpha ;\,\cdot\,)$ is continuous, it is
uniformly continuous on $\overline \Omega \times\overline{
B}(0,\lambda)$, and we define the modulus of continuity
$\omega_{\eta,\lambda} : Q' \times \mathbb R^+ \rightarrow \mathbb
R^+$ by

\begin{eqnarray*}&\omega_{\eta,\lambda}(y_\alpha,t):= \sup
\Big\{\,|W^{\eta,\lambda}(x,y_\alpha;\xi)
-W^{\eta,\lambda}(x',y_\alpha;\xi')|,~\text{where} \\
&\hspace{5cm} (x,\xi), \, (x',\xi') \in \overline \Omega \times
\overline{ B}(0,\lambda) \text{and}~~|(x;\xi)-(x';\xi')| \leq t
\Big\}.\end{eqnarray*}

 Then

$$\left\{
\begin{array}{l}
\omega_{\eta,\lambda}(\cdot,t) \text{ is lower semicontinuous for
all }t \in \mathbb
R^+,\\
\\
\omega_{\eta,\lambda}(y_\alpha,\cdot) \text{ is continuous and
increasing for all }y_\alpha
\in Q',\\\\
\omega_{\eta,\lambda}(y_\alpha,0)=0 \text{ for all }y_\alpha \in
Q',\end{array} \right.$$

\noindent and

\begin{equation}\label{modcont}
|W^{\eta,\lambda}(x,y_\alpha;\xi)-W^{\eta,\lambda}(x',y_\alpha;\xi')|
\leq \omega_{\eta,\lambda}(y_\alpha,|x-x'|+|\xi-\xi'|)
\end{equation}

\noindent for all  $(x,\xi),\; (x',\xi') \in \overline \Omega
\times \overline{ B}(0,\lambda).$ The first property is a
consequence of the fact that the supremum of continuous functions
is lower semicontinuous, while the other ones are classical
properties of moduli of continuity.

For all $t \in \mathbb R^+$, we extend
$\omega_{\eta,\lambda}(\cdot,t)$ to $\mathbb R^2$ by
$Q'$-periodicity. Since $W^{\eta,\lambda}(x,\,\cdot\,;\xi)$ is
$Q'$-periodic, inequality  (\ref{modcont}) holds for all $y_\alpha
\in \mathbb R^2$\!. Consequently, for every $(x_\alpha,x_3) \in
[Q'_{i,h} \times I] \cap R^\lambda_n $ and every $x_\alpha' \in
Q'_{i,h}$,

\begin{eqnarray*}
&&\left|W^{\eta,\lambda}\left(x_\alpha,x_3,\frac{x_\alpha}{\varepsilon_n};\overline
\xi +D_\alpha u_n(x) \Big| \frac{1}{\varepsilon_n} D_3
u_n(x)\right)\right.\\
&&\hspace{3.0cm}\left.
-W^{\eta,\lambda}\left(x_\alpha',x_3,\frac{x_\alpha}{\varepsilon_n};\overline
\xi +D_\alpha u_n(x) \Big| \frac{1}{\varepsilon_n} D_3
u_n(x)\right)\right|\\
&&\hspace{0.5cm} \leq
\omega_{\eta,\lambda}\left(\frac{x_\alpha}{\varepsilon_n},|x_\alpha-x_\alpha
'|\right)\\&& \hspace{0.5cm}\leq
\omega_{\eta,\lambda}\left(\frac{x_\alpha}{\varepsilon_n},\frac{\sqrt
2 \delta}{h}\right).
\end{eqnarray*}

After integration  in $(x_\alpha,x_3,x_\alpha')$ and summation, we
get

\begin{eqnarray*}
&&\sum_{i=1}^{h^2} \frac{h^2}{\delta^2} \int_{Q'_{i,h}}
\left\{\int_{R^\lambda_n \cap [ Q'_{i,h} \times
I]}\left|W^{\eta,\lambda}\left(x_\alpha,x_3,\frac{x_\alpha}{\varepsilon_n};\overline
\xi +D_\alpha u_n(x) \Big| \frac{1}{\varepsilon_n} D_3
u_n(x)\right)\right.\right.\\
&&\hspace{3.0cm} \left.\left.
-W^{\eta,\lambda}\left(x_\alpha',x_3,\frac{x_\alpha}{\varepsilon_n};\overline
\xi +D_\alpha u_n(x) \Big| \frac{1}{\varepsilon_n} D_3
u_n(x)\right)\right| dx \right\} dx_\alpha'\\
&& \hspace{0.5cm}\leq  2
\int_{Q'(x_\alpha^0,\delta)}\omega_{\eta,\lambda}\left(\frac{x_\alpha}{\varepsilon_n},
\frac{\sqrt 2 \delta}{h}\right)dx_\alpha.
\end{eqnarray*}

Riemann-Lebesgue's Lemma applied to the $Q'$-periodic function
$\omega_{\eta,\lambda}(\,\cdot\,,\sqrt 2 \delta/h)$ yields,

$$\lim_{n\to +\infty}2
\int_{Q'(x_\alpha^0,\delta)}\omega_{\eta,\lambda}\left(\frac{x_\alpha}{\varepsilon_n},
\frac{\sqrt 2 \delta}{h}\right)dx_\alpha = 2\delta ^2
\int_{Q'}\omega_{\eta,\lambda}\left(x_\alpha,\frac{\sqrt 2
\delta}{h}\right)dx_\alpha,$$

and by Beppo-Levi's Monotone Convergence Theorem

$$\lim_{h \to +\infty}2\delta ^2
\int_{Q'}\omega_{\eta,\lambda}\left(x_\alpha,\frac{\sqrt 2
\delta}{h}\right)dx_\alpha=0.$$

Hence

\begin{eqnarray*}
&\ds \mathcal I_{\{\varepsilon_n\}}(v;Q'(x_\alpha^0,\delta)) \geq \\
&\ds \limsup_{\lambda,\eta,h,n}\sum_{i=1}^{h^2}
\frac{h^2}{\delta^2} \!\!\int_{ Q'_{i,h}} \!\! \left\{ \int_{[
Q'_{i,h} \times I] \cap R^\lambda_n
}\!\!\!\!W^{\eta,\lambda}\!\left(\!x_\alpha',x_3,\frac{x_\alpha}{\varepsilon_n};
\overline{\xi} +D_\alpha u_n \Big| \frac{1}{\varepsilon_n} D_3
u_n\!\right)dx \right\} dx_\alpha'.
\end{eqnarray*}

Define the following sets which depend on all parameters
$(\eta,\lambda,i,h,n)$ :\\

$$\left\{
\begin{array}{l}
T := \{ (x_\alpha',x_\alpha,x_3) \in Q'_{i,h}\times Q'_{i,h}
\times I : (x_\alpha',x_3)
\in K_\eta \text{ and } (x_\alpha,x_3) \in R^\lambda_n\},\\
\\
T_1 := \{ (x_\alpha',x_\alpha,x_3) \in Q'_{i,h}\times Q'_{i,h}
\times I :
(x_\alpha',x_3) \not \in K_\eta \text{ and } (x_\alpha,x_3) \in R^\lambda_n\},\\
\\
T_2:= \{ (x_\alpha',x_\alpha,x_3) \in Q'_{i,h}\times Q'_{i,h}
\times I : (x_\alpha,x_3) \not \in R^\lambda_n\},
\end{array}
\right.$$

\noindent and note that $Q'_{i,h} \times Q'_{i,h} \times I= T \cup
T_1 \cup T_2$. Since $W(\cdot,y_\alpha;\,\cdot\,)$ and
$W^{\eta,\lambda}(\cdot,y_\alpha;\,\cdot\,)$ coincide on
$K_\eta\times \overline{B}(0,\lambda)$, we have

\begin{eqnarray}
&&\mathcal I_{\{\varepsilon_n\}}(v;Q'(x_\alpha^0,\delta))\nonumber\\
&&\hspace{0.5cm} \geq  \limsup_{\lambda,\eta,h,n} \sum_{i=1}^{h^2}
\frac{h^2}{\delta^2} \int_{T} W^{\eta,\lambda}
\left(x_\alpha',x_3,\frac{x_\alpha}{\varepsilon_n};\overline \xi
+D_\alpha u_n(x) \Big| \frac{1}{\varepsilon_n} D_3
u_n(x)\right)dx \, dx_\alpha' \nonumber\\
&&\hspace{1.0cm} = \limsup_{\lambda,\eta,h,n} \sum_{i=1}^{h^2}
\frac{h^2}{\delta^2} \int_{T}
W\left(x_\alpha',x_3,\frac{x_\alpha}{\varepsilon_n};\overline \xi
+D_\alpha u_n(x) \Big| \frac{1}{\varepsilon_n} D_3 u_n(x)\right)dx
\, dx_\alpha'.\label{E}
\end{eqnarray}

We will prove that the corresponding terms over $T_1$ and $T_2$
are zero. Indeed, in view of (\ref{mes}) and the $p$-growth
condition $(H_{3})$,

\begin{eqnarray}
&&\sum_{i=1}^{h^2}  \frac{h^2}{\delta^2}\int_{T_1}
W\left(x_\alpha',x_3,\frac{x_\alpha}{\varepsilon_n};\overline \xi
+D_\alpha u_n(x) \Big| \frac{1}{\varepsilon_n} D_3
u_n(x)\right)dx\,
dx_\alpha'\nonumber\\
&& \hspace{3.0cm}  \leq \quad \sum_{i=1}^{h^2}
\frac{h^2}{\delta^2} \, \mathcal L^2(Q'_{i,h}) \, \mathcal
L^3([Q'_{i,h} \times I] \setminus K_\eta)\, \beta
(1+\lambda^p)\nonumber\\
&& \hspace{3.0cm}<  \quad \beta (1+\lambda^p)\eta
\xrightarrow[\eta \to 0]{} 0,\label{E1}
\end{eqnarray}

\noindent uniformly in $(n,h)$. The bound from above in $(H_{3})$,
the equi-integrability of $\left\{\left|\left(D_\alpha u_n \big|
\frac{1}{\varepsilon_n} D_3 u_n\right)\right|^p\right\}$ and
(\ref{lambda}) imply that

\begin{eqnarray}
&& \sum_{i=1}^{h^2}  \frac{h^2}{\delta^2}\int_{T_2}
W\left(x_\alpha',x_3,\frac{x_\alpha} {\varepsilon_n};\overline \xi
+D_\alpha u_n(x) \Big| \frac{1}{\varepsilon_n}
D_3 u_n(x)\right)dx\, dx_\alpha' \nonumber\\
&& \hspace{2.0cm}\leq \quad \sum_{i=1}^{h^2}  \frac{h^2}
{\delta^2} \, \mathcal L^2(Q'_{i,h}) \, \beta \int_{[Q'_{i,h}
\times I] \setminus R^\lambda_n}\left(1+\left|\left(D_\alpha
u_n(x) \Big| \frac{1}{\varepsilon_n} D_3
u_n(x)\right)\right|^p\right)dx\nonumber\\
&& \hspace{2.0cm} = \quad \beta \int_{[Q'(x_\alpha^0,\delta)
\times I] \setminus R^\lambda_n}\left(1+\left|\left(D_\alpha
u_n(x) \Big| \frac{1}{\varepsilon_n} D_3
u_n(x)\right)\right|^p\right)dx \xrightarrow[\lambda \to
+\infty]{} 0,\label{E2}
\end{eqnarray}

uniformly in $(\eta,n,h)$. Thus, in view of (\ref{E}), (\ref{E1}),
(\ref{E2}), Fatou's Lemma yields

\begin{eqnarray*}
& \ds \mathcal I_{\{\varepsilon_n\}}(v;Q'(x_\alpha^0,\delta))\\
& \ds \geq   \! \limsup_{h \to +\infty} \limsup_{n \to +\infty}
\sum_{i=1}^{h^2} \frac{h^2}{\delta^2}\int_{Q'_{i,h}} \!\!
\int_{Q'_{i,h} \times
I}\!\!W\left(x_\alpha',x_3,\frac{x_\alpha}{\varepsilon_n};\overline
\xi +D_\alpha u_n \Big| \frac{1}{\varepsilon_n} D_3
u_n\right)\!dx \,dx_\alpha'\\
& \ds \geq  \!\limsup_{h \to +\infty} \sum_{i=1}^{h^2}
\frac{h^2}{\delta^2} \!\! \int_{Q'_{i,h}}  \liminf_{n \to
+\infty}\int_{Q'_{i,h} \times
I}\!\!W\left(x_\alpha',x_3,\frac{x_\alpha}{\varepsilon_n};\overline
\xi +D_\alpha u_n \Big| \frac{1}{\varepsilon_n} D_3
u_n\right)\!dx\, dx_\alpha'.
\end{eqnarray*}

Fix $x_\alpha' \in Q'_{i,h}$ such that $W_{\rm
hom}(x_\alpha';\overline \xi)$ is well defined and set
$Z(x;\xi):=W(x_\alpha',x_3,x_\alpha;\xi)$. It is easy to check
that $Z$ is a Carath\'eodory integrand hence, applying  Theorem
4.2 of  Braides, Fonseca and Francfort \cite{Bra&Fo&Fr}, we get
since $u_n \to 0$ in $L^p(Q'(x_\alpha^0,\delta) \times I;\Rb^3)$,

$$2\frac{\delta^2}{h^2} \overline Z(\overline \xi) \leq \liminf_{n \to +\infty}
\int_{Q'(x_\alpha^0,\delta)\times
I}Z\left(\frac{x_\alpha}{\varepsilon_n},x_3;\overline \xi
+D_\alpha u_n(x)\Big| \frac{1}{\varepsilon}D_3u_n(x)\right)\,
dx,$$ where

\begin{eqnarray*}
\overline Z(\overline \xi) & := & \inf_{T>0,\;
\varphi}\Big\{\int_{(0,T)^2 \times I}
Z(x;\overline \xi +D_\alpha \varphi(x)|D_3 \varphi(x))\,dx : \\
& &  \hspace{2.0cm} \varphi \in W^{1,p}((0,T)^2 \times I;\Rb^3),
\quad \varphi=0\text{ on }\partial (0,T)^2 \times I\Big\}.
\end{eqnarray*}

In view of the previous formula together with (\ref{whom}) and
Remark \ref{infwhom}, we have that $\overline Z(\overline
\xi)=W_{\rm hom}(x_\alpha';\overline \xi)$. Then

$$\liminf_{n \to +\infty}\int_{Q'_{i,h} \times
I}W\left(x_\alpha',x_3,\frac{x_\alpha}{\varepsilon_n};\overline
\xi +D_\alpha u_n(x) \Big| \frac{1}{\varepsilon_n} D_3
u_n(x)\right)dx \geq \frac{2 \delta^2}{h^2}\;  W_{\rm
hom}(x_\alpha';\overline \xi),$$

\noindent and so

\begin{eqnarray*}\mathcal I_{\{\varepsilon_n\}}(v;Q'(x_\alpha^0,\delta))
&\geq& \limsup_{h \to +\infty}\sum_{i=1}^{h^2}
\frac{h^2}{\delta^2} \int_{Q'_{i,h}} \frac{2 \delta^2}{h^2}\;
W_{\rm hom}(x_\alpha';\overline \xi) dx_\alpha'\\
& = &2 \int_{Q'(x_\alpha^0,\delta)}W_{\rm hom}(x_\alpha';\overline
\xi) dx_\alpha'.\end{eqnarray*}

Dividing both sides of the previous inequality by $\delta^2$ and
passing to the limit when $\delta \searrow 0^+$, we obtain by
(\ref{C1)}) and (\ref{lpg})

$$W_{\{\varepsilon_n\}}(x_\alpha^0;\overline \xi) \geq   W_{\rm hom}
(x_\alpha^0;\overline \xi).$$ \hfill $\blacksquare$

\vspace{0.2cm} \vspace{0.5 cm}

\begin{proposition}\label{prop-ult}

$W_{\{\e_n\}}(x_{\alpha};\overline\xi)= W_{\rm
hom}(x_{\alpha};\overline\xi)$ a.e.\! $x_{\alpha}\in \o$ and all
$\overline\xi\in \Rb^{3\times 2}.$
\end{proposition}

{\it Proof.}  Let $E$ be the intersection of the set $L$ (see
Definition \ref{setL}) with the subset of points $x_{\alpha}^0 \in
\o$ where $W_{\{\e_n\}}(x_{\alpha}^{0};\,\cdot\,)$ and $W_{\rm
hom}(x_{\alpha}^{0};\,\cdot\,)$ are continuous (see Lemma
\ref{cont}). Then $\mathcal L^2(\o \setminus E)=0$ and in view of
Lemma \ref{ineq1} and \ref{ineq2}, we have that for all
$x_\alpha^0 \in E$ and for all $\overline \xi \in \mathbb Q^{3
\times 2}$,

$W_{\{\varepsilon_n\}}(x_\alpha^0;\overline \xi) = W_{\rm hom}
(x_\alpha^0;\overline \xi)$. Since
$W_{\{\varepsilon_n\}}(x_\alpha^0;\, \cdot\, )$ and $W_{\rm
hom}(x_\alpha^0;\, \cdot \,)$ are continuous for all $x_\alpha^0
\in E$, the equality $W_{\{\varepsilon_n\}}(x_\alpha^0;\overline
\xi) = W_{\rm hom} (x_\alpha^0;\overline \xi)$ holds true for
$x_\alpha^0 \in E$ and all $\overline \xi \in \mathbb R^{3 \times
2}$\!.

\hfill $\blacksquare$

\begin{coro}\label{coro-ult}

$\Gamma(L^p(A\times I)){\text-}\debaixodolim {}{\e} \mathcal
I_{\varepsilon}(\,\cdot\,;A) =\mathcal I_{\rm hom}(\,\cdot\,;A)$
for all $A\in {\cal A}(\o),$ where  $\mathcal I_{\rm
hom}(\cdot;A)$ is the functional defined in (\ref{claim}).
\end{coro}

{\it Proof.} From Proposition \ref{prop-ult}  we can conclude that
$\mathcal I_{\rm hom}(\cdot;A)$ is well defined and

 $$\Gamma(L^p(A\times I)){\text-}\debaixodolim {}{n}
\mathcal I_{\varepsilon_n}(\,\cdot\,;A) =\mathcal I_{\rm
hom}(\,\cdot\,;A)$$

\noindent  for all $A\in {\cal A}(\o)$ (see Remark \ref{dem2}).
Since this limit does not depend upon the extracted subsequence,
in view of Proposition 7.11 in Braides and Defranceschi
\cite{Bra&De}, the whole sequence $\{\mathcal
I_\varepsilon(\,\cdot\,;A)\}_{\e>0}$ $\Gamma(L^p(A\times
I))$-converges to $\mathcal I_{\rm hom}(\,\cdot\,;A)$ for each
$A\in {\cal A}(\o)$. \hfill
$\blacksquare$\\

The proof of Theorem \ref{jf-m} comes as a consequence of
Corollary \ref{coro-ult} taking $A=\o.$

\section{Appendix}

\noindent We now prove a technical result on extension of
Carath\'eodory functions that was useful in the proof of Lemma
\ref{ineq2}. The argument used is very close to that of Theorem 1,
Section 1.2 in Evans and Gariepy \cite{Ev&Ga}.

\begin{lemma}\label{ext}

Let $\Omega \subset \Rb^N$ be a bounded open set and $f : \Omega
\times \Rb^m \times \Rb^{d \times N} \to \Rb$ a function such that
$$\left\{\begin{array}{l}
f(x,\,\cdot\,;\,\cdot\,) \text{ is continuous for a.e. }x \in \Omega;\\\\
f( \,\cdot\, y ;\xi) \text{ is }\mathcal L^N
\text{-measurable for all } y\in \Rb^{N}\,\,{\rm and }\,\, \xi \in \Rb^N;\\\\
f(x,\,\cdot\,;\xi) \text{ is }(0,1)^m\text{-periodic for a.e.\! }x
\in \Omega \text{ and all }\xi \in \Rb^{d \times N}.
\end{array}\right.$$

Assume also that there exists $\beta >0$ and $1\leq p<\infty$ such
that

$$\frac{1}{\beta}|\xi|^p-\beta \leq f(x,y;\xi) \leq
\beta(1+|\xi|^p),\quad\text{ for a.e. \! } x \in \Omega \text{ and
all }(y,\xi) \in \Rb^m \times \Rb^{d \times N}.$$

Then for any $\eta>0$ and $\lambda>0$ there exist a compact set
$K_\eta \subset \Omega$ and a function $f^{\eta,\lambda} : \mathbb
R^N \times \mathbb R^m \times \mathbb R^{d \times N} \to \Rb$ such
that

$$\left\{\begin{array}{l}\mathcal L^N(\Omega \setminus K_\eta) < \eta,\\\\
f^{\eta,\lambda}(x,y;\xi)=f(x,y;\xi) \text{ for all }(x,y,\xi) \in
K_\eta
\times \mathbb R^m \times\overline{B}(0,\lambda),\\\\
f^{\eta,\lambda}(\, \cdot\,,y;\,\cdot\,) \text{ is continuous for
all }y \in\mathbb R^m,\\\\
f^{\eta,\lambda}(x,\,\cdot\,;\xi) \text{ is continuous and
}(0,1)^m\text{-periodic for all }(x,\xi) \in \mathbb R^N \times
\mathbb R^{d \times N},
\end{array}\right.$$

and

\begin{equation}\label{ContExt}
-\beta \leq f^{\eta,\lambda}(x,y;\xi) \leq
\beta(1+\lambda^p),\quad \text{ for all }(x,y,\xi) \in \mathbb R^N
\times \mathbb R^m \times \mathbb R^{d \times N}.
\end{equation}
\end{lemma}

\noindent {\it Proof. } Since $f$ is a Carath\'eodory function, by
Scorza Dragoni's Theorem (see Ekeland and Teman \cite{Ek&Te}) for
all $\eta>0$ there exists a compact set $K_\eta \subset \Omega$
satisfying $\mathcal L^N(\Omega \setminus K_\eta)<\eta$ and such
that $f$ is continuous on $K_\eta \times \mathbb R^m \times
\mathbb R^{d \times N}$. Let $C^{\eta,\lambda}:= K_\eta \times
\overline{B}(0,\lambda)\equiv C$ (to simplify notation) and $
U^{\eta,\lambda}:=(\Rb^N \times \Rb^{d \times N}) \setminus
C^{\eta,\lambda}\equiv U$. Fix $(s,\gamma) \in C$, and for all
$(x,\xi) \in U$ set

$$ u^{\eta,\lambda}_{(s,\gamma)}(x,\xi):=\max\left\{ 2-\frac{|(s,\gamma)-(x,\xi)|}{{\rm
dist} ((x,\xi),C)} ,0\right\}\equiv u_{(s,\gamma)}(x,\xi).$$

\noindent Clearly

$$\left\{\begin{array}{l}
u_{(s,\gamma)} \text{ is continuous on }U,\\\\
0\leq u_{(s,\gamma)} \leq 1,\\\\
u_{(s,\gamma)}(x,\xi)=0 \text{ if and only if }
|(s,\gamma)-(x,\xi)| \geq2{\rm dist} ((x,\xi),C).
\end{array}\right.$$

Let $\{s^{\eta}_j\}_{j \geq 1}\equiv \{s_j\}_{j \geq 1}$ and
$\{\gamma^{\lambda}_j\}_{j \geq 1}\equiv \{\gamma_j\}_{j \geq 1}$
be  countable dense families  in  $K_\eta$ and
$\overline{B}(0,\lambda)$, respectively. Define

$$\sigma^{\eta,\lambda}(x,\xi):=\sum_{j \geq 1}2^{-j}u_{(s_j,\gamma_j)}(x,\xi)\equiv
\sigma(x,\xi) \quad \text{ for all}\,\,(x,\xi) \in U.$$

Since $\sigma$ is the uniform limit of a sequence of continuous
functions in $U$, then $\sigma$ is continuous in $U$. Moreover,
for all $(x,\xi) \in U$ $$0< \sigma(x,\xi) \leq 1.$$

 Indeed,assume that $\sigma(x,\xi)=0$ for some $(x,\xi) \in U$. Then, for
all $j \geq 1$, $u_{(s_j,\gamma_j)}(x,\xi)=0$ and thus
$|(s_j,\gamma_j)-(x,\xi)| \geq 2\, {\rm dist} ((x,\xi),C)$. The
density of $\{s_j,\gamma_j\}$ in $C$ yields that

$$|(s,\gamma)-(x,\xi)| \geq 2\, {\rm dist} ((x,\xi),C)$$
\noindent for all $(s,\gamma) \in C$. We obtain a contradiction if
we choose $(s,\gamma)$ to be those points of $C$ such that ${\rm
dist} ((x,\xi),C)= {\rm dist}((x,\xi),(s,\gamma))$ so
$\sigma(x,\xi)>0$ for all $(x,\xi)\in U$. Consequently, the
function

$$(x,\xi) \mapsto v_k(x,\xi) \equiv v_k^{\eta,\lambda}(x,\xi):=
\frac{2^{-k}u_{(s_k,\gamma_k)}(x,\xi)}{\sigma(x,\xi)}$$ is well
defined and continuous in $U$. Moreover it satisfies that

$$0 \leq v_k(x,\xi) \leq 1,\quad \sum_{k \geq 1}v_k(x,\xi)=1 \quad \text{ for
all}\,\, (x,\xi)\in U.$$ Fix $y \in \mathbb R^m$ and define the
continuous extension of $f(\cdot,y;\,\cdot\,)$ outside $C$ as
$$f^{\eta,\lambda}(x,y;\xi)=\left\{\begin{array}{ll}
\ds f(x,y,\xi) & \text{ if }(x,\xi) \in C,\\
\ds \sum_{k \geq 1} v_k(x,\xi)\, f(s_k,y;\gamma_k) & \text{ if
}(x,\xi)\in U.
\end{array}\right.$$

Obviously, we have $f^{\eta,\lambda}(x,y;\xi)=f(x,y;\xi)$ for all
$(x,y,\xi) \in K_\eta \times \mathbb R^m \times
\overline{B}(0,\lambda)$. On the other hand, if $(x,y,\xi)$ is
such that $(x,\xi)\in U$, in view of the $p$-growth and the
$p$-coercivity condition on $f$ we get  $$-\beta \leq
f^{\eta,\lambda}(x,y;\xi) \leq \sum_{k \geq 1}v_k(x,\xi)
\beta(1+|\gamma_k|^p) \leq \beta(1+\lambda^p).$$ Since we have

\begin{equation}\label{613}
\sup_{y \in\Rb^m, \,(x,\xi)\in U} \Big[ \sum_{k\geq
n}\left|2^{-k}u_{(s_k,\gamma_k)}(x,\xi)f(s_k,y;\gamma_k)\right|
\Big] \leq \beta(1+\lambda^p)\sum_{k\geq n} 2^{-k}\xrightarrow[n
\to +\infty]{}0,
\end{equation}

\noindent then the function

$$(x,y,\xi) \mapsto \sum_{k\geq 1} 2^{-k}u_{(s_k,\gamma_k)}(x,\xi)f(s_k,y;\gamma_k)$$

\noindent is continuous on $\{(x,y,\xi): (x,\xi) \in U, \; y \in
\Rb^m\}$. In particular, for all $(x,\xi) \in \Rb^N \times \Rb^{d
\times N}$ the function $f^{\eta,\lambda}(x,\,\cdot\,;\xi)$ is
continuous. Furthermore, $f^{\eta,\lambda}(x,\,\cdot\,;\xi)$  is
$(0,1)^m$-periodic because if $\mathbf i \in \mathbb Z^m$ then for
$(x,\xi)\in U$

$$f^{\eta,\lambda}(x,y+\mathbf i ;\xi)=\sum_{k \geq 1} v_k(x,\xi)\,
f(s_k,y+\mathbf i;\gamma_k)=\sum_{k \geq 1} v_k(x,\xi)\,
f(s_k,y;\gamma_k)=f^{\eta,\lambda}(x,y;\xi).$$

Finally we prove the continuity of $f^{\eta,\lambda}(\cdot,
y_;\cdot)$. By (\ref{613}) it suffices to show that for all $(a,A)
\in C$

$$\lim_{U \ni (x,\xi) \to (a,A)}f^{\eta,\lambda}(x,y;\xi)=f(a,y;A).$$

\noindent As $\{(s_j,\gamma_j)\}_{j \geq 1}$ is dense in $C$ and
$f(\cdot,y;\cdot)$ is continuous on $C$, for every $\e>0$ there
exists $\delta>0$ such that $|f(a,y;A) - f(s_j,y;\gamma_j)| < \e$
for all $j \geq 1$ with $|(a,A)-(s_j,\gamma_j)| <\delta$. Assume
that $|(x,\xi)-(a,A)| < \delta /4$ and suppose that $j \geq 1$ is
such that $|(a,A)-(s_j,\gamma_j)| \geq \delta$. Then

$$\delta \leq |(a,A)-(s_j,\gamma_j)| \leq |(a,A)-(x,\xi)| + |(x,\xi)-(s_j,\gamma_j)|\leq
\frac{\delta}{4}+|(x,\xi)-(s_j,\gamma_j)|,$$

\noindent and thus

$$|(x,\xi)-(s_j,\gamma_j)| \geq \frac{3\delta}{4} > 2 |(a,A)-(x,\xi)| \geq 2\, {\rm
dist}((x,\xi),C).$$

\noindent Consequently, $v_j(x,\xi)=0$ if $j$ is such that
$|(a,A)-(s_j,\gamma_j)| \geq \delta$, and so

$$|f^{\eta,\lambda}(x,y;\xi)-f(a,y;A)| \leq \sum_{j \geq 1, \, |(a,A)-(s_j,\gamma_j)|
< \delta} v_j(x,\xi)| f(s_j,y;\gamma_j)-f(a,y;A)| <\e,$$

\noindent because non zero terms of the sum are those which
satisfy $|f(a,y;A) - f(s_j,y;\gamma_j)| < \e$. The continuity of
$f^{\eta,\lambda}(\cdot,y;\cdot)$ now follows.
\hfill $\blacksquare$\\

\noindent {\bf Acknowledgements:}  The authors would like to thank
Irene Fonseca and Gilles Francfort  for  suggesting them this
research
 work, and for their fruitful comments and suggestions.  We wish to thank  the Center
for
 Nonlinear Analysis (NSF Grant No. 0405343) for
its hospitality and support. The research of M. Ba\'{\i}a  was
partially supported by Funda\c{c}\~{a}o para a Ci\^{e}ncia e
Tecnologia (Grant PRAXIS XXI SFRH
$\hspace{-0.1cm}\backslash\hspace{-0.1cm}$ BD
$\hspace{-0.1cm}\backslash\hspace{-0.1cm}$ 1174
$\hspace{-0.1cm}\backslash\hspace{-0.1cm}$ 2000), Fundo Social
Europeu, the Department of Mathematical Sciences of Carnegie
Mellon University and its Center for Nonlinear Analysis.

\vspace{0.5cm}

\begin{center}

\begin{small}

Jean-Fran\c{c}ois Babadjian\\

 \textsc{L.P.M.T.M., Universit\'e Paris Nord, 93430, Villetaneuse, France}\\

\textit{E-mail address}: {\bf jfb@galilee.univ-paris13.fr}

\vspace{0.5cm}

Margarida Ba\'{\i}a\\

 \textsc{Departamento de Matem\'{a}tica, Instituto Superior T\'{e}cnico,}\\

 \textsc{1049-001 Lisboa, Portugal}\\

\textit{E-mail address}: {\bf mbaia@math.ist.utl.pt}

\end{small}

\end{center}

\end{document}